\documentclass[11pt]{amsart}
\usepackage{amsmath}
\usepackage{amsxtra}
\usepackage{amscd}
\usepackage{amsthm}
\usepackage{amsfonts}
\usepackage{amssymb}
\usepackage{eucal}

\usepackage{latexsym}

\newtheorem{cor}{Corollary}
\newtheorem{lem}{Lemma}
\newtheorem{prop}{Proposition}

\newtheorem{thm}{Theorem}

\theoremstyle{remark}
\newtheorem{remark}{Remark}

\theoremstyle{definition}

\numberwithin{equation}{section}

\newcommand{\thmref}[1]{Theorem~\ref{#1}}

\newcommand{\secref}[1]{Sect.~\ref{#1}}
\newcommand{\lemref}[1]{Lemma~\ref{#1}}
\newcommand{\propref}[1]{Proposition~\ref{#1}}

\newcommand{\nc}{\newcommand}
\nc{\renc}{\renewcommand}
\nc{\ssec}{\subsection}
\nc{\sssec}{\subsubsection}
\nc{\on}{\operatorname}

\nc\ol{\overline}
\nc\ul{\underline}
\nc\wt{\widetilde}
\nc\tboxtimes{\wt{\boxtimes}}
\nc{\wh}{\widehat}
\nc{\mc}{\mathcal}

\nc{\CM}{{\mathcal M}}
\nc{\CN}{{\mathcal N}}
\nc{\CF}{{\mathcal F}}
\nc{\D}{{\mathcal D}}
\nc{\CQ}{{\mathcal Q}}
\nc{\CY}{{\mathcal Y}}
\nc{\CX}{{\mathcal X}}
\nc{\CG}{{\mathcal G}}
\nc{\CE}{{\mathcal E}}
\nc{\CC}{{\mathcal C}}
\nc{\CO}{{\mathcal O}}
\renc{\CC}{{\mathcal C}}
\nc{\CT}{{\mathcal T}}
\nc{\CK}{{\mathcal K}}
\nc{\CS}{{\mathcal S}}
\nc{\CH}{{\mathcal H}}
\nc{\CU}{{\mathcal U}}
\nc{\CV}{{\mathcal V}}
\nc{\CA}{{\mathcal A}}
\nc{\CB}{{\mathcal B}}
\nc{\CW}{{\mathcal W}}
\nc{\CL}{{\mathcal L}}
\nc{\CP}{{\mathcal P}}
\nc{\CI}{{\mathcal I}}
\nc{\CJ}{{\mathcal J}}
\nc{\CR}{{\mathcal R}}

\nc{\BA}{{\mathbb{A}}}
\nc{\BC}{{\mathbb{C}}}
\nc{\BG}{{\mathbb{G}}}
\nc{\BM}{{\mathbb{M}}}
\nc{\BN}{{\mathbb{N}}}
\nc{\BP}{{\mathbb{P}}}
\nc{\BR}{{\mathbb{R}}}
\nc{\BZ}{{\mathbb{Z}}}
\nc{\BV}{{\mathbb{V}}}
\nc{\BW}{{\mathbb{W}}}
\nc{\BS}{{\mathbb{S}}}
\nc{\BD}{{\mathbb{D}}}
\nc{\BQ}{{\mathbb{Q}}}
\nc{\BL}{{\mathbb{L}}}
\renc{\BW}{{\mathbb{W}}}

\nc{\fa}{{\mathfrak{a}}}
\nc{\fb}{{\mathfrak{b}}}
\nc{\fg}{{\mathfrak{g}}}
\nc{\fgl}{{\mathfrak{gl}}}
\nc{\fh}{{\mathfrak{h}}}
\nc{\fj}{{\mathfrak{j}}}
\nc{\fm}{{\mathfrak{m}}}
\nc{\fl}{{\mathfrak{l}}}
\nc{\fn}{{\mathfrak{n}}}
\nc{\fu}{{\mathfrak{u}}}
\nc{\fp}{{\mathfrak{p}}}
\nc{\fr}{{\mathfrak{r}}}
\nc{\fs}{{\mathfrak{s}}}
\nc{\fsl}{{\mathfrak{sl}}}

\nc{\hsl}{{\widehat{\mathfrak{sl}}}}
\nc{\hgl}{{\widehat{\mathfrak{gl}}}}
\nc{\hg}{{\widehat{\mathfrak{g}}}}
\nc{\hb}{{\widehat{\mathfrak{b}}}}
\nc{\hn}{{\widehat{\mathfrak{n}}}}

\nc{\fA}{{\mathfrak{A}}}
\nc{\fB}{{\mathfrak{B}}}
\nc{\fO}{{\mathfrak{O}}}
\nc{\fD}{{\mathfrak{D}}}
\nc{\fE}{{\mathfrak{E}}}
\nc{\fF}{{\mathfrak{F}}}
\nc{\fG}{{\mathfrak{G}}}
\nc{\fK}{{\mathfrak{K}}}
\nc{\fL}{{\mathfrak{L}}}
\nc{\fC}{{\mathfrak{C}}}
\nc{\fM}{{\mathfrak{M}}}
\nc{\fN}{{\mathfrak{N}}}
\nc{\fH}{{\mathfrak{H}}}
\nc{\fP}{{\mathfrak{P}}}
\nc{\fU}{{\mathfrak{U}}}
\nc{\fV}{{\mathfrak{V}}}
\nc{\fZ}{{\mathfrak{Z}}}
\nc{\fz}{{\mathfrak{z}}}

\nc{\bc}{{\mathbf{c}}}
\nc{\bd}{{\mathbf{d}}}
\nc{\bh}{{\mathbf{h}}}
\nc{\be}{{\mathbf{e}}}
\nc{\bj}{{\mathbf{j}}}
\nc{\bn}{{\mathbf{n}}}
\nc{\bp}{{\mathbf{p}}}
\nc{\bg}{{\mathbf{g}}}
\nc{\bq}{{\mathbf{q}}}
\nc{\bs}{{\mathbf{s}}}
\nc{\bu}{{\mathbf{u}}}
\nc{\bv}{{\mathbf{v}}}
\nc{\bx}{{\mathbf{x}}}
\nc{\by}{{\mathbf{y}}}
\nc{\bw}{{\mathbf{w}}}
\nc{\bA}{{\mathbf{A}}}
\nc{\bK}{{\mathbf{K}}}
\nc{\bB}{{\mathbf{B}}}
\nc{\bC}{{\mathbf{C}}}
\nc{\bD}{{\mathbf{D}}}
\nc{\bH}{{\mathbf{H}}}
\nc{\bM}{{\mathbf{M}}}
\nc{\bN}{{\mathbf{N}}}
\nc{\bV}{{\mathbf{V}}}
\nc{\bW}{{\mathbf{W}}}
\nc{\bL}{{\mathbf{L}}}
\nc{\bU}{{\mathbf{U}}}
\nc{\bX}{{\mathbf{X}}}
\nc{\bI}{{\mathbf{I}}}
\nc{\bZ}{{\mathbf{Z}}}
\nc{\bS}{{\mathbf{S}}}

\nc{\sA}{{\mathsf{A}}}
\nc{\sB}{{\mathsf{B}}}
\nc{\sC}{{\mathsf{C}}}
\nc{\sD}{{\mathsf{D}}}
\nc{\sF}{{\mathsf{F}}}
\nc{\sH}{{\mathsf{H}}}
\nc{\sG}{{\mathsf{G}}}
\nc{\sK}{{\mathsf{K}}}
\nc{\sM}{{\mathsf{M}}}
\nc{\sO}{{\mathsf{O}}}
\nc{\sQ}{{\mathsf{Q}}}
\nc{\sP}{{\mathsf{P}}}
\nc{\sV}{{\mathsf{V}}}
\nc{\sZ}{{\mathsf{Z}}}
\nc{\sfp}{{\mathsf{p}}}
\nc{\sr}{{\mathsf{r}}}
\nc{\sg}{{\mathsf{g}}}
\nc{\sk}{{\mathsf{k}}}
\nc{\ssf}{{\mathsf{f}}}
\nc{\ssh}{{\mathsf{h}}}
\nc{\sse}{{\mathsf{e}}}
\nc{\sfb}{{\mathsf{b}}}
\nc{\sfc}{{\mathsf{c}}}
\nc{\sd}{{\mathsf{d}}}

\nc{\Av}{\on{Av}}
\nc{\act}{\on{act}}
\nc{\Hom}{\on{Hom}}
\nc{\End}{\on{End}}
\nc{\Lie}{\on{Lie}}
\nc{\Loc}{\on{Loc}}
\nc{\IC}{\on{IC}}
\nc{\Aut}{\on{Aut}}
\nc{\rk}{\on{rk}}
\nc{\Sh}{\on{Sh}}
\nc{\Perv}{\on{Perv}}
\nc{\pos}{{\on{pos}}}
\nc{\Conv}{\on{Conv}}
\nc{\Sph}{\on{Sph}}
\nc{\Sym}{\on{Sym}}
\nc{\Rep}{{\mc R}ep(\cG)}
\nc{\RepH}{{\mc R}ep(H)}
\nc{\Fun}{\on{Fun}}
\nc{\Id}{\on{Id}}
\nc{\id}{\on{id}}
\renc{\mod}{\on{-mod}}

\nc{\oG}{\overset{\circ}{G}{}}
\nc{\oGB}{{\overset{\circ}{G/B}{}}}
\nc{\oGN}{{\overset{\circ}{G/N}{}}}
\nc{\uBC}{\underline{\BC}}

\nc{\crit}{{\on{crit}}}
\nc{\reg}{{\on{reg}}}
\nc{\nilp}{{\on{nilp}}}
\nc{\ord}{\on{ord}}
\nc{\nil}{\wt{\on{reg}}}
\nc{\mb}{\mathbf}
\nc{\ren}{\on{ren}}
\nc{\res}{\on{res}}
\nc{\RS}{{\on{RS}}}
\nc{\Dist}{\on{Dist}}
\nc{\semiinf}{{\frac{\infty}{2}}}
\nc{\semiinfi}{{\frac{\infty}{2}+i}}
\nc{\semiinfb}{{\frac{\infty}{2}+\bullet}}
\nc{\torsemiinf}{{\overset{\semiinf}\otimes}}
\nc{\Hitch}{\on{Hitch}}

\nc{\hl}{\overset{\leftarrow}h}
\nc{\hr}{\overset{\rightarrow}h}
\nc\Dh{\widehat{\D}}
\nc{\Gr}{\on{Gr}}
\nc{\Fl}{\on{Fl}}
\nc{\Flt}{\wt{\Fl}{}}
\nc{\Pic}{\on{Pic}}
\nc{\Bun}{\on{Bun}}

\nc{\bDR}{\mathbf {DR}}
\nc{\uV}{\underline{V}}
\nc{\arrowtimes}{\overset{\to}\otimes}
\nc{\hattimes}{\widehat\otimes}
\nc{\larrowtimes}{\overset{\leftarrow}\otimes}
\nc{\shriektimes}{\overset{!}\otimes}
\nc{\startimes}{\overset{*}\otimes}
\nc{\sCliff}{\mathsf {Cliff}}
\nc{\sSpin}{\mathsf {Spin}}

\nc{\one}{{\mathbf{1}}}

\nc\Spec{\on{Spec}}
\nc{\Pro}{\on{Pro}}
\nc{\QCoh}{\on{QCoh}}
\nc{\uHom}{\underline{\on{Hom}}}
\nc{\RHom}{\on{RHom}}
\nc{\uRHom}{\underline{\on{RHom}}}
\nc{\CHom}{{\mathcal Hom}}
\nc{\uCHom}{\underline{{\mathcal Hom}}}
\nc{\uCRHom}{\underline{{\mathcal R}{\mathcal Hom}}}

\nc{\cg}{\check \fg}
\nc{\Op}{\on{Op}_{\cg}}
\nc{\nOp}{\on{Op}^{\nilp}_{\cg}}
\nc{\nMOp}{\on{MOp}^{\nilp}_{\cg}}
\nc{\rOp}{\on{Op}^{\reg}_{\cg}}

\nc{\tg}{\wt{\check \fg}}
\nc{\cn}{\check \fn}
\nc{\tn}{\wt{\cn}}
\nc{\cG}{\check G}
\nc{\cB}{\check B}
\nc{\cb}{\check \fb}
\nc{\MOp}{\on{MOp}_{\check \cg}}
\nc{\cN}{\CN_{\cG}}
\nc{\tN}{\wt{\CN}_{\cG}}
\nc{\dIsom}{{\mathsf{Isom}}_{\Op}}
\nc{\disom}{{\mathsf{isom}}_{\Op}}
\nc{\Kdv}{{\mathsf{Isom}}_{\Op^\reg}}
\nc{\kdv}{{\mathsf{isom}}_{\Op^\reg}}
\nc{\Isom}{{\mathsf{Isom}}_{\on{Op}_\fg}}
\nc{\isom}{{\mathsf{isom}}_{\on{Op}_\fg}}

\nc{\wcosta}{j_{\wt{w},*}}
\nc{\wsta}{j_{\wt{w},!}}
\nc{\wcost}{j_{w,*}}
\nc{\wst}{j_{w,!}}

\nc{\epsi}{{\mathbf e}^\psi}
\nc{\epsip}{{\mathbf e}^{\psi'}}

\nc{\Ppi}{{\mathbf \Pi}}

\nc{\hCO}{{\hat{\CO}}}
\nc{\hCK}{{\hat{\CK}}}

\nc{\CPreg}{\CP_{G,\on{Op}^\reg}}
\nc{\CPBreg}{\CP_{B,\on{Op}^\reg}}
\nc{\CPnilp}{\CP_{G,\on{Op}^\nilp}}
\nc{\CPBnilp}{\CP_{B,\on{Op}^\nilp}}
\nc{\CPla}{\CP_{G,\on{Op}_{\cla}}}
\nc{\CPBla}{\CP_{B,\on{Op}_{\cla}}}

\nc{\Cat}{\hg_\crit\mod^{I,m}_\nilp}
\nc{\Catf}{{}^f\hg_\crit\mod^{I,m}_\nilp}
\nc{\DCat}{D^b(\hg_\crit\mod_\nilp)^{I^0}}
\nc{\DCatf}{{}^f D^b(\hg_\crit\mod_\nilp)^{I^0}}
\nc{\Catr}{\hg_\crit\mod^{I,m}_\reg}
\nc{\Catrf}{{}^f\hg_\crit\mod^{I,m}_\reg}
\nc{\DCatr}{D^b(\hg_\crit\mod_\reg)^{I^0}}
\nc{\DCatrf}{{}^f D^b(\hg_\crit\mod_\reg)^{I^0}}

\nc{\ch}{\mbox{ch}}
\nc{\Z}{{\mathbb Z}}
\nc{\C}{{\mathbb C}}
\nc{\pone}{{\mathbb C}{\mathbb P}^1}
\nc{\pa}{\partial}
\nc{\F}{{\mathcal F}}
\nc{\arr}{\rightarrow}
\nc{\larr}{\longrightarrow}
\nc{\al}{\alpha}
\nc{\ri}{\rangle}
\nc{\lef}{\langle}
\nc{\W}{{\mathcal W}}
\nc{\la}{\lambda}
\nc{\ep}{\epsilon}
\nc{\su}{\widehat{{\mathfrak s}{\mathfrak l}}_2}
\nc{\sw}{{\mathfrak s}{\mathfrak l}}
\nc{\g}{{\mathfrak g}}
\nc{\h}{{\mathfrak h}}
\nc{\n}{{\mathfrak n}}
\nc{\N}{\widehat{\n}}
\nc{\De}{\Delta}
\nc{\gt}{\widetilde{\g}}
\nc{\Ga}{\Gamma}
\nc{\z}{{\mathfrak Z}}
\nc{\La}{\Lambda}
\nc{\cri}{_{\kappa_c}}
\nc{\kk}{h^\vee}
\nc{\sun}{\widehat{\sw}_N}
\nc{\si}{\sigma}
\nc{\el}{\ell}
\nc{\bi}{\bibitem}
\nc{\om}{\omega}
\nc{\ds}{\displaystyle}
\nc{\dzz}{\frac{dz}{z}}
\nc{\Res}{\on{Res}}
\nc{\Cal}{\mathcal}
\nc{\bb}{{\mathfrak b}}
\nc{\ot}{\otimes}
\nc{\R}{{\mc R}}
\nc{\yy}{{\mc Y}}
\nc{\ga}{\gamma}

\nc{\us}{\underset}
\nc{\opl}{\oplus}
\nc{\beq}{\begin{equation}}
\nc{\Fq}{{\mathbb F}_q}
\nc{\Mq}{{\mathcal M}}
\nc{\lan}{\langle}
\nc{\ran}{\rangle}

\nc{\Vect}{\on{Vect}}
\nc{\ghat}{\wh\fg}
\nc{\T}{\mc T}
\nc{\Tloc}{\T^\g_{\on{loc}}}
\nc{\vac}{|0\ran}
\nc{\Wick}{{\mb :}}
\nc{\delz}{\partial_z}
\nc{\K}{{\cali K}}
\nc{\cali}{\mathcal}
\nc{\li}{\mathfrak l}
\nc{\lt}{\widetilde{\li}}
\nc{\astar}{a^*}
\nc{\cA}{{\mc A}}
\nc{\ka}{\kappa}

\nc{\OO}{{\mc O}}
\nc{\AutO}{\on{Aut}\OO}
\nc{\DerO}{\on{Der}\OO}
\nc{\DerpO}{\on{Der}_+\OO}
\nc{\Au}{{\mc A}ut}
\nc{\mf}{\mathfrak}
\nc{\V}{{\mc V}}
\nc{\hh}{\wh{\h}}

\nc{\pp}{{\mathfrak p}}
\nc{\mm}{{\mathfrak m}}
\nc{\rr}{{\mathfrak r}}
\nc{\ket}{\rangle}
\nc{\zz}{\fZ^{\reg}}
\nc{\gr}{\on{gr}}
\nc{\Spe}{\on{Spec}}
\nc{\rv}{\crho}
\nc{\can}{\on{can}}
\nc{\Db}{{\mathbb D}}
\nc{\ww}{w}

\nc{\RR}{\on{R}}
\nc{\PPi}{{\mathbf \Pi}}
\nc{\M}{{\mathbb M}}
\nc{\Mv}{{\mathbb M}^\vee}
\nc{\VV}{{\mathbb V}}
\nc{\bsl}{\backslash}

\nc{\bchi}{{\mathbf {\chi}}}
\nc{\anch}{{\mathbf {anch}}}

\nc{\cla}{{\check{\la}}}
\nc{\cmu}{{\check{\mu}}}
\nc{\crho}{{\check{\rho}}}
\nc{\com}{{\check{\omega}}}
\nc{\DD}{{\mc D}}
\nc{\E}{{\mc E}}
\nc{\Ll}{{\mc L}}

\nc{\ConnX}{\on{Conn}_H(\omega_X^{\crho})}
\nc{\ConHX}{\on{Conn}_{\check{H}}(\omega_X^{\rho})}
\nc{\ConnD}{\on{Conn}_H(\omega_{\D}^{\crho})}
\nc{\ConHD}{\on{Conn}_{\check{H}}(\omega_{\D}^{\rho})}
\nc{\ConnDt}{\on{Conn}_H(\omega_{\D^\times}^{\crho})}
\nc{\ConHDt}{\on{Conn}_{\check{H}}(\omega_{\D^\times}^{\rho})}

\nc{\Ind}{\on{Ind}}

\nc{\CTop}{{\mathcal Top}}

\nc{\ppart}{(\!(t)\!)}

\nc{\qu}{/\!/}

\nc{\gen}{\on{gen}}

\nc{\chal}{\check\al}
\nc{\zmod}{\zz\mod}
\nc{\gmod}{\hg_\crit\mod}
\nc{\zlamod}{\zz_\la\mod}

\begin{document}

\title{Weyl modules and opers without monodromy}

\author{Edward Frenkel}\thanks{The research of E.F. was supported by
the DARPA Program ``Focus Areas in Fundamental Mathematics''.}

\address{Department of Mathematics, University of California,
  Berkeley, CA 94720, USA}

\email{frenkel@math.berkeley.edu}

\author{Dennis Gaitsgory}

\address{Department of Mathematics, Harvard University,
Cambridge, MA 02138, USA}

\email{gaitsgde@math.harvard.edu}

\date{June 2007}

\begin{abstract}
We prove that the algebra of endomorphisms of a Weyl module of
critical level is isomorphic to the algebra of functions on the space
of monodromy-free opers on the disc with regular singularity and
residue determined by the highest weight of the Weyl module. This
result may be used to test the local geometric Langlands
correspondence proposed in our earlier work.
\end{abstract}

\maketitle

\section{Introduction}

Let $\fg$ be a simple finite-dimensional Lie algebra. For an invariant
inner product $\kappa$ on $\fg$ (which is unique up to a scalar) define
the central extension $\hg_\kappa$ of the formal loop algebra $\fg
\otimes \BC\ppart$ which fits into the short exact sequence
$$
0 \to \BC {\mb 1} \to \hg_\kappa \to \fg \otimes \BC\ppart \to 0.
$$
This sequence is split as a vector space, and the commutation
relations read
\begin{equation}    \label{KM rel}
[x \otimes f(t),y \otimes g(t)] = [x,y] \otimes f(t) g(t) - (\kappa(x,y)
\on{Res} f dg) {\mb 1},
\end{equation}
and ${\mb 1}$ is a central element. The Lie algebra $\hg_\kappa$ is
the {\em affine Kac--Moody algebra} associated to $\kappa$. We will
denote by $\hg_\kappa\mod$ the category of {\em discrete}
representations of $\hg_{\kappa}$ (i.e., such that any vector is
annihilated by $\fg \otimes t^n\BC[[t]]$ for sufficiently large $n$),
on which ${\mb 1}$ acts as the identity.

Let $U_{\kappa}(\hg)$ be the quotient of the universal enveloping
algebra $U(\hg_{\kappa})$ of $\hg_{\kappa}$ by the ideal generated by
$({\mb 1}-1)$. Define its completion $\wt{U}_{\kappa}(\hg)$ as follows:
$$
\wt{U}_{\kappa}(\hg) = \underset{\longleftarrow}\lim \;
U_{\kappa}(\hg)/U_{\kappa}(\hg) \cdot (\fg \otimes t^n\BC[[t]]).
$$
It is clear that $\wt{U}_{\kappa}(\hg)$ is a topological algebra,
whose discrete continuous representations are the same as objects of
$\hg_\kappa\mod$.

Let $\kappa_{\crit}$ be the {\em critical} inner product on $\fg$
defined by the formula
$$
\kappa_{\crit}(x,y) = - \frac{1}{2} \on{Tr} (\on{ad} (x) \circ \on{ad}
(y)).
$$
In what follows we will use the subscript ``$\crit$'' instead of
$\kappa_\crit$.

Let $\check G$ be the group of adjoint type whose Lie algebra
$\check \fg$ is Langlands dual to $\fg$ (i.e., the Cartan matrix of
$\check\fg$ is the transpose of that of $\fg$).

Let $\fZ_{\g}$ be the center of $\wt{U}_{\crit}(\hg)$. According
to a theorem of \cite{FF,F:wak}, $\fZ_{\g}$ is isomorphic to the
algebra $\on{Fun} \on{Op}_{\check \fg}(\D^\times)$ of functions on the
space of $\check \fg$-opers on the punctured disc $\D^\times =
\on{Spec}(\C\ppart)$ (see \cite{BD,FG:local} and \secref{first
section} for the definition of opers).

It is interesting to understand how $\fZ_{\g}$ acts on various
$\hg_\crit$-modules. The standard modules are the Verma modules and
the Weyl modules. They are obtained by applying the induction functor
\begin{align*}
\on{Ind}: \fg\mod &\to \hg_\crit\mod, \\
M &\mapsto U(\hg_\crit) \underset{U(\fg[[t]] \oplus {\mb 1})}\otimes M,
\end{align*}
where $\fg[[t]]$ acts on $M$ via the projection $\fg[[t]] \to \fg$ and
${\mb 1}$ acts as the identity.

For $\la \in \h^*$ let $M_\la$ be the Verma module over $\fg$ with
highest weight $\la$. The corresponding $\hg_\crit$-module $\BM_\la =
\on{Ind}(M_\la)$ is the {\em Verma module} of critical level with
highest weight $\la$.

For a dominant integral weight $\la$ let $V_\la$ be the irreducible
finite-dimensional $\g$-module with highest weight $\la$. The
corresponding $\hg_\crit$-module $\BV_\la = \on{Ind}(V_\la)$ is the
{\em Weyl module} of critical level with highest weight $\la$. The
module $\BV_0 = \Ind(\BC_0)$ is also called the vacuum module.

It was proved in \cite{FF,F:wak} that the algebra of
$\hg_\crit$-endomorphisms of $\BV_0$ is isomorphic to the algebra
$\on{Fun}\Op^\reg$ of functions on the space $\Op^\reg$ of $\check
\fg$-opers on the disc $\D = \on{Spec}(\BC[[t]])$. Moreover, there is
a commutative diagram
$$
\CD
\fZ_{\g} @>{\sim}>> \Fun\on{Op}_{\check \fg}(\D^\times) \\
@VVV   @VVV  \\
\on{End}_{\hg_\crit}(\BV_0) @>{\sim}>> \Fun\Op^\reg
\endCD
$$

We have shown in \cite{FG:local}, Corollary 13.3.2, that a similar
result holds for the Verma modules as well: the algebra of
$\hg_\crit$-endomorphisms of $\BM_\la$ is isomorphic to
$\on{Fun}\Op^{\RS,\varpi(-\lambda-\rho)}$, where
$\Op^{\RS,\varpi(-\lambda-\rho)}$ is the space of $\cg$-opers
on $\D^\times$ with regular singularity and residue
$\varpi(-\lambda-\rho)$, $\varpi$ being the natural projection $\h^*
\to \on{Spec} \on{Fun}(\h^*)^W$ (see \cite{FG:local}, Sect. 2.4, for a
precise definition). In addition, there is an analogue of the above
commutative diagram for Verma modules.

In this paper we consider the Weyl modules $\BV_\la$. In
\cite{FG:local}, Sect. 2.9, we defined the subspace $\Op^{\la,\reg}
\subset \Op^{\RS,\varpi(-\lambda-\rho)}$ of $\la$-{\em regular opers}
(we recall this definition below). Its points are those opers in
$\Op^{\RS,\varpi(-\lambda-\rho)}$ which have trivial monodromy and are
therefore $\cG(\hCK)$ gauge equivalent to the trivial local system on
$\D^\times$. In particular, $\Op^{0,\reg} = \Op^{\reg}$. According to
\lemref{no monodromy} below, the disjoint union of $\Op^{\la,\reg}$,
where $\la$ runs over the set $P^+$ of dominant integral weights of
$\g$, is precisely the locus of $\cg$-opers on $\D^\times$ with
trivial monodromy. The main result of this paper is the following
theorem, which generalizes the description of $\on{End}_{\ghat_\crit}
\BV_0$ from \cite{FF,F:wak} to the case of an arbitrary dominant
integral weight $\la$.

\begin{thm}    \label{weyl end}
For any dominant integral weight $\la$ the center $\fZ_{\g}$ maps
surjectively onto $\on{End}_{\ghat_\crit} \BV_\la$, and we have the
following commutative diagram
\begin{equation}    \label{weyl diag}
\begin{CD}
\fZ_{\g} @>{\sim}>> \on{Fun}\on{Op}_{\cG}(D^\times)
\\ @VVV @VVV \\ \on{End}_{\ghat_{\crit}} \BV_\la
@>{\sim}>> \on{Fun}\Op^{\la,\reg}
\end{CD}
\end{equation}
\end{thm}

For $\g=\sw_2$ this follows from Prop. 1 of \cite{F:icmp}. This
statement was also independently conjectured by A. Beilinson and
V. Drinfeld (unpublished).

In addition, we prove that $\BV_\la$ is a free module over
$\on{End}_{\ghat_\crit} \BV_\la \simeq \on{Fun}\Op^{\la,\reg}$.

\thmref{weyl end} has important consequences for the local geometric
Langlands correspondence proposed in \cite{FG:local}. According to our
proposal, to each ``local Langlands parameter'' $\sigma$, which is a
$\cG$--local system on the punctured disc $\D^\times$ (or
equivalently, a $\cG$-bundle with a connection on $\D^\times$), there
should correspond a category $\CC_\sigma$ equipped with an action of
$G\ppart$.

Now let $\chi$ be a fixed $\cg$-oper on $\D^\times$, which we regard
as a character of the center $\fZ_{\g}$. Consider the full
subcategory $\hg_\crit\mod_\chi$ of the category $\hg_\crit\mod$ whose
objects are $\ghat_\crit$-modules, on which the $\fZ_{\g}$ acts
according to this character. This category carries a canonical action
of the ind-group $G\ppart$ via its adjoint action on $\hg_\crit$. We
proposed in \cite{FG:local} that $\hg_\crit\mod_\chi$ should be
equivalent to the sought-after category $\CC_\sigma$, where $\sigma$
is the $\cG$-local system underlying the oper $\chi$. This entails a
far-reaching corollary that the categories $\hg_\crit\mod_{\chi_1}$
and $\hg_\crit\mod_{\chi_2}$ for two different opers
$\chi_1$ and $\chi_2$ are equivalent if the underlying local systems
of $\chi_1$ and $\chi_2$  are isomorphic to each other.

In particular, consider the simplest case when $\sigma$ is the trivial
local system. Then, by \lemref{no monodromy}, $\chi$ must be a point
of $\Op^{\la,\reg}$ for some $\la \in P^+$. \thmref{weyl end} implies
that the quotient $\BV_\la(\chi)$ of the Weyl module $\BV_\la$ by the
central character corresponding to $\chi$ is non-zero. Therefore
$\BV_\la(\chi)$ is a non-trivial object of $\hg_\crit\mod_\chi$ and
also of the corresponding $G[[t]]$-equivariant category
$\hg_\crit\mod_\chi^{G[[t]]}$. In the case when $\la=0$ we have proved
in \cite{FG:exact} that this is a unique irreducible object of
$\hg_\crit\mod_\chi^{G[[t]]}$ and that this category is in fact
equivalent to the category of vector spaces. Therefore we expect the
same to be true for all other values of $\la$. This will be proved in
a follow-up paper.

\medskip

The paper is organized as follows. In \secref{first section} we recall
the relevant notions of opers, Cartan connections and Miura
transformation. In \secref{proof main} we explain the strategy of the
proof of the main result, \thmref{weyl end}, and reduce it to two
statements, \thmref{factors} and \propref{injective}. We then prove
\propref{injective} assuming \thmref{factors} in
\secref{exactness}. Our argument is based on the exactness of the
functor of quantum Drinfeld--Sokolov reduction, which we derive from
\cite{Ar}. In \secref{char z} we compute the characters of the algebra
of functions on $\Op^{\la,\reg}$ and of the semi-infinite cohomology
of $\BV_\la$ (they turn out to be the same). We then give two
different proofs of \thmref{factors} in \secref{proof}. This completes
the proof of the main result. We also show that the natural map from
the Weyl module $\BV_\la$ to the corresponding Wakimoto module is
injective and that $\BV_\la$ is a free module over its endomorphism
algebra.

\section{Some results on opers}    \label{first section}

In this section we recall the relevant notions of opers, Cartan
connections and Miura transformation, following
\cite{BD,F:wak,FG:local}, where we refer the reader for more details.

%We will also compute the character of the algebra of
%functions on the space $\on{Op}_{\g}^{\cla,\reg}$ which will be used
%in the proof of our main theorem.

Let $\g$ be a simple Lie algebra and $G$ the corresponding algebraic
group of adjoint type. Let $B$ be its Borel subgroup and $N = [B,B]$
its unipotent radical, with the corresponding Lie algebras $\fn
\subset \fb \subset \fg$.

Let $X$ be a smooth curve, or the disc $\D = \on{Spec} (\hCO)$, where
$\hCO$ is a one-dimensional smooth complete local ring, or the
punctured disc $\D^\times=\Spec (\hCK)$, where $\hCK$ is the field of
fractions of $\hCO$.

Following Beilinson and Drinfeld (see \cite{BD}, Sect. 3.1, and
\cite{BD1}), one defines a $\fg$-{\em oper} on $X$ to be a triple
$(\CF_G,\nabla,\CF_B)$, where $\CF_G$ is a principal $G$-bundle
$\CF_G$ on $X$, $\nabla$ is a connection on $\CF_G$, and $\CF_B$ is a
$B$-reduction of $\CF_G$ which is transversal to $\nabla$, in the
sense explained in the above references and in \cite{FG:local},
Sect. 1.1. We note that the transversality condition allows us to
identify canonically the $B$-bundles $\F_B$ underlying all opers.

More concretely, opers on the punctured disc $\D^\times$ may be
described as follows. Let us choose a trivialization of $\CF_B$ and a
coordinate $t$ on the disc $\D$ such that $\hCO = \C[[t]]$ and $\hCK =
\C\ppart$. Let us choose a nilpotent subalgebra $\fn_-$, which is in
generic position with $\fb$ and a set of simple root generators
$f_\imath, \imath \in I$, of $\fn_-$. Then a $\g$-oper on $\D^\times$
is, by definition, an equivalence class of operators of the form
\begin{equation}    \label{oper on disc}
\nabla = \pa_t + \sum_{\imath \in I} f_\imath + \bv, \qquad \bv \in
\fb(\hCK),
\end{equation}
with respect to the action of the group $N(\hCK)$ by gauge
transformations. It is known that this action is free and the
resulting set of equivalence classes is in bijection with
$\hCK^{\oplus \ell}$, where $\ell = \on{rank}(\g)$. Opers may be
defined in this way over any base, and this allows us to define the
ind-affine scheme $\on{Op}_{\g}(\D^\times)$ of $\g$-opers on
$\D^\times$ (it is isomorphic to an inductive limit of affine
spaces).

Let $\check{P}^+$ be the set of dominant integral coweights of
$\g$. For $\cla \in \check{P}^+$ we define a $\fg$-{\em oper with
$\cla$-nilpotent singularity} as an equivalence class of operators
\begin{equation}    \label{cla opers}
\nabla = \pa_t + \sum_{\imath \in I}
t^{\langle\alpha_\imath,\cla\rangle} \cdot f_\imath + \bv(t) +
\frac{v}{t}, \qquad \bv(t) \in \fb(\hCO), v \in \fn,
\end{equation}
with respect to the action of the group $N(\hCO)$ by gauge
transformations (see \cite{FG:local}, Sect. 2.9). The corresponding
scheme is denoted by $\on{Op}_{\g}^{\cla,\nilp}$. According to Theorem
2.9.1 of \cite{FG:local}, the natural map $\on{Op}_{\g}^{\cla,\nilp}
\to \on{Op}_{\g}(\D^\times)$ is injective, and its image is equal to
the space of $\g$-opers with regular singularity and residue
$\varpi(-\cla-\crho)$ (here $\crho$ is the half-sum of positive
coroots of $\g$ and $\varpi$ is the projection $\h \to \on{Spec}
\on{Fun}(\h)^W$).

Now we define the space $\on{Op}_{\g}^{\cla,\reg}$ of $\cla$-{\em
regular opers} as the subscheme of $\on{Op}_{\g}^{\cla,\nilp}$
corresponding to those operators \eqref{cla opers} which satisfy $v =
0$ (so that $\nabla$ is regular at $t=0$). In particular, if $\cla=0$,
then $\on{Op}_{\g}^{0,\reg}$ is just the space of regular $\cg$-opers
on the disc $\D$. The geometric significance of $\cla$-opers is
explained by the following

\begin{lem}    \label{no monodromy}
Suppose that a $\fg$-oper $\chi = (\CF_G,\nabla,\CF_B)$ on $\D^\times$
is such that the corresponding $G$-local system is trivial (in other
words, the corresponding operator \eqref{oper on disc} is in the
$G(\hCK)$ gauge equivalence class of $\nabla_0 = \pa_t$). Then $\chi
\in \on{Op}_{\g}^{\cla,\reg}$ for some $\cla \in \check{P}^+$.
\end{lem}

\begin{proof}
It is clear from the definition that any oper in
$\on{Op}_{\g}^{\cla,\reg}$ is regular on the disc $\D$ and is
therefore $G(\hCK)$ gauge equivalent to the trivial connection
$\nabla_0 = \pa_t$.

Now suppose that we have an oper $\chi = (\CF_G,\nabla,\CF_B)$ on
$\D^\times$ such that the corresponding $G$-local system is trivial.
Then $\nabla$ is $G(\hCK)$ gauge equivalent to a regular connection on
$\D$. We have the decomposition $G(\hCK) = G(\hCO) B(\hCK)$. The gauge
action of $G(\hCO)$ preserves the space of regular connections (in
fact, it acts transitively on it). Therefore if an oper connection
$\nabla$ is gauge equivalent to a regular connection under the action
of $G(\hCK)$, then its $B(\hCK)$ gauge equivalence class must contain
a regular connection. The oper condition then implies that this gauge
class contains a connection operator of the form \eqref{cla opers}
with $\bv(0) = 0$, for some dominant integral coweight
$\cla$. Therefore $\chi \in \on{Op}_{\g}^{\cla,\reg}$.
\end{proof}

Let us choose a coordinate $t$ on $\D$. The vector field $L_0 = - t
\pa_t$ then acts naturally on $\on{Op}_{\g}^{\cla,\reg}$ and defines a
$\Z$-grading on the algebra of functions on it. In \secref{char z} we
will compute the character of the algebra $\on{Fun}
\on{Op}_{\g}^{\cla,\reg}$ of functions on $\on{Op}_{\g}^{\cla,\reg}$
with respect to this grading.

\medskip

Next, we introduce the space of $H$-connections and the Miura
transformation. 

Let $X$ be as above. Denote by $\omega_X$ the $\C^\times$-torsor
corresponding to the canonical line bundle on $X$. Let
$\omega^{\crho}_X$ be the push-forward of $\omega_X$ to an $H$-torsor
via the homomorphism $\crho:\C^\times\to H$. We denote by $\ConnX$ the
affine space of all connections on $\omega^{\crho}_X$. In particular,
$\ConnDt$ is an inductive limit of affine spaces $\ConnD^{\ord_k}$ of
connections with pole of order $\leq k$. We will use the notation
$\ConnD^\RS$ for $\ConnD^{\ord_1}$. A connection $\ol\nabla \in
\ConnD^\RS$ has a well-defined residue, which is an element of
$\h$. For $\cmu \in \h$ we denote by $\ConnD^{\RS,\cmu}$ the subspace
of $\ConnD^\RS$ consisting of connections with residue $\cmu$.

The {\em Miura transformation} is a morphism 
\begin{equation}    \label{MT}
\on{MT}: \ConnDt \to \on{Op}_{\g}(\D^\times),
\end{equation}
introduced in \cite{DS} (see also \cite{F:wak} and \cite{FG:local},
Sect. 3.3). It can be described as follows. If we choose a coordinate
$t$ on $\D$, then we trivialize $\omega_{\D}$ and hence
$\omega^{\crho}_{\D}$. A point of $\ConnDt$ is then represented by an
operator
$$
\ol\nabla = \pa_t + {\mb u}(t), \qquad {\mb u}(t) \in \h(\hCK).
$$
We associate to $\ol\nabla$ the $\fg$-oper which is the $N(\hCK)$
gauge equivalence class of the operator
$$
\nabla = \pa_t + \sum_{\imath \in I} f_{\imath} + {\mb u}(t).
$$

The following result is a corollary of Proposition 3.5.4 of
\cite{FG:local}. Let $\F_{B,0}$ be the fiber at $0 \in \D$ of the
$B$-bundle $\F_B$ underlying all $\fg$-opers. We denote by
$N_{\F_{B,0}}$ the $\F_{B,0}$-twist of $N$.

\begin{lem}    \label{prin N bdle}
Let $\cla$ be a dominant integral coweight of $\g$. The image of
$\ConnD^{\RS,-\cla}$ in $\on{Op}_{\g}(\D^\times)$ under the Miura
transformation is equal to $\on{Op}_{\g}^{\cla,\reg}$. Moreover, the
map $\ConnD^{\RS,-\cla} \to \on{Op}_{\g}^{\cla,\reg}$ is a principal
$N_{\F_{B,0}}$-bundle over $\on{Op}_{\g}^{\cla,\reg}$.
\end{lem}

In particular, this implies that the scheme $\on{Op}_{\g}^{\cla,\reg}$
is smooth and in fact isomorphic to an infinite-dimensional
(pro)affine space.

\section{Proof of the main theorem}    \label{proof main}

%\subsection{Outline of the proof}

Our strategy of the proof of \thmref{weyl end} will be as follows: we
will first construct natural maps
\begin{equation}    \label{hom of Ends}
\on{End}_{\ghat_\crit} \BM_\la \to \on{End}_{\ghat_\crit} \BV_\la \to 
H^{\frac{\infty}{2}}(\fn_+\ppart,\fn_+[[t]],\BV_\la \otimes
\Psi_0).
\end{equation}
We already know from \cite{FG:local} that
$$
\on{End}_{\ghat_\crit} \BM_\la \simeq \on{Fun} \Op^{\la,\nilp}.
$$
We will show that the corresponding composition
$$
\on{Fun} \Op^{\la,\nilp} \to 
H^{\frac{\infty}{2}}(\fn_+\ppart,\fn_+[[t]],\BV_\la \otimes
\Psi_0)
$$
factors as follows:
$$
\on{Fun} \Op^{\la,\nilp} \twoheadrightarrow \on{Fun} \Op^{\la,\reg}
\overset{\sim}\longrightarrow
H^{\frac{\infty}{2}}(\fn_+\ppart,\fn_+[[t]],\BV_\la \otimes \Psi_0),
$$
and that the map
$$
\BV_\la \to H^{\frac{\infty}{2}}(\fn_+\ppart,\fn_+[[t]],\BV_\la
\otimes \Psi_0)
$$
is injective. This will imply \thmref{weyl end}.

As a byproduct, we will obtain an isomorphism
$$
H^{\frac{\infty}{2}}(\fn_+\ppart,\fn_+[[t]],\BV_\la \otimes
\Psi_0) \simeq \on{End}_{\ghat_\crit} \BV_\la
$$
and find that the first map in \eqref{hom of Ends} is surjective.

\subsection{Homomorphisms of $\ghat_\crit$-modules}    \label{hom
  ghat}

Let us now proceed with the proof and construct the maps
\eqref{hom of Ends}.

Note that a $\ghat_\crit$-endomorphism of $\BM_\la$ is uniquely
determined by the image of the highest weight vector, which must be a
vector in $\BM_\la$ of weight $\la$ annihilated by the Lie subalgebra
$$
\wh\n_+ = (\n_+ \otimes 1) \oplus (\g \otimes t\C[[t]]).
$$
This is the Lie algebra of the prounipotent proalgebraic group $I^0 =
[I,I]$, where $I$ is the Iwahori subgroup of $G\ppart$.
For a $\ghat_\crit$-module $M$ we denote the space of such vectors
by $M^{\wh\n_+}_\la$.

According Corollary 13.3.2 of \cite{FG:local}, we have
$$
\on{End}_{\ghat_\crit} \BM_\la = (\BM_\la)^{\wh\n_+}_\la \simeq
\on{Fun} \Op^{\la,\nilp} = \on{Fun}
\on{Op}_{\cG}^{\on{RS},\varpi(-\la-\rho)}
$$
(see \secref{first section}).

Likewise, any endomorphism of $\BV_\la$ is uniquely determined by the
image of the generating subspace $V_\la$. This subspace therefore
defines a $\g[[t]]$-invariant vector in $(\BV_\la \otimes
V_\la^*)^{\g[[t]]}$. Note that for any $\g$-integrable module $M$ we
have an isomorphism
$$
(M \otimes V_\la^*)^{\g[[t]]} \simeq M^{\wh\n_+}_\la.
$$
Therefore we have
\begin{equation}    \label{end isomo}
\on{End}_{\ghat_\crit} \BV_\la = (\BV_\la \otimes
V_\la^*)^{\g[[t]]} = (\BV_\la)^{\wh\n_+}_\la.
\end{equation}

The canonical surjective homomorphism
$$
\BM_\la \twoheadrightarrow \BV_\la
$$
of $\ghat_\crit$-modules gives rise to a map $(\BM_\la)^{\wh\n_+}_\la
\to (\BV_\la)^{\wh\n_+}_\la$. We obtain the following commutative
diagram:
\begin{equation}    \label{important diagram1}
\begin{CD}
\fZ_{\g} @>>> \on{End}_{\ghat_\crit} \BM_\la
@>>> \on{End}_{\ghat_\crit} \BV_\la \\ @VVV @VV{\sim}V
@VV{\sim}V \\ \on{Fun} \Op(\D^\times) @>>> \on{Fun} \Op^{\la,\nilp}
@>>> ?
\end{CD}
\end{equation}

\subsection{The functor of semi-infinite cohomology}    \label{cdo}

Define the character
\begin{equation}    \label{Psi0}
\Psi_0: \n_+\ppart \to \C
\end{equation}
by the formula
$$
\Psi_0(e_{\al,n}) = \begin{cases} 1, & \on{if} \al = \al_\imath, n=-1,
  \\ 0, & \on{otherwise}. \end{cases}
$$
We have the functor of semi-infinite cohomology (the $+$ quantum
Drinfeld--Sokolov reduction) from the category of
$\ghat_\crit$-modules to the category of graded vector spaces,
\begin{equation}    \label{n+}
M \mapsto H^{\frac{\infty}{2}+\bullet}(\fn_+\ppart,\fn_+[[t]],M
\otimes \Psi_0),
\end{equation}
introduced in \cite{FF,FKW} (see also \cite{vertex}, Ch. 15, and
\cite{FG:local}, Sect. 18; we follow the notation of the latter).

Let $M$ be a $\ghat_\crit$-module. Consider the space
$M^{\wh\n_+}_\la$ of $\wh\n_+$-invariant vectors in $M$ of highest
weight $\la$.

\begin{lem}    \label{map}
We have functorial maps
\begin{equation}   \label{si coh inj}
M^{\wh\n_+}_\la \to
H^{\frac{\infty}{2}}(\fn_+\ppart,\fn_+[[t]],M \otimes \Psi_0).
\end{equation}
\end{lem}

\begin{proof}
Consider the complex $C^{\bullet}(M)$ computing the above
semi-infinite cohomology (see, e.g., \cite{vertex}, Sect. 15.2). It
follows from the definition that the standard Chevalley complex
computing the cohomology of the Lie algebra $\n_+[[t]]$ with
coefficients in $M$ embeds into $C^{\bullet}(M)$. Therefore we obtain
functorial maps
$$
M^{\wh\n_+}_\la \to M^{\n_+[[t]]} \to
H^{\frac{\infty}{2}}(\fn_+\ppart,\fn_+[[t]],\BV_\la \otimes \Psi_0).
$$
\end{proof}

Introduce the notation
\begin{align*}
\fz^{\la,\nilp} &= \on{Fun} \Op^{\la,\nilp}, \\
\fz^{\la,\reg} &= \on{Fun} \Op^{\la,\reg}
\end{align*}

Our goal is to prove that
$$
(\BV_\la)^{\wh\n_+}_\la \simeq \fz^{\la,\reg}.
$$
The proof will be based on analyzing the composition
\begin{equation}    \label{comp}
\fz^{\la,\nilp} \to (\BV_\la)^{\wh\n_+}_\la \to
H^{\frac{\infty}{2}}(\fn_+\ppart,\fn_+[[t]],\BV_\la \otimes \Psi_0),
\end{equation}
where the first map is obtained from the diagram \eqref{important
diagram1}, and the second map from \lemref{map}. We will use the
following two results.

\begin{thm}    \label{factors}
%There is a canonical isomorphism
%$$
%H^{\frac{\infty}{2}}(\fn_+\ppart,\fn_+[[t]],\BV_\la \otimes \Psi_0)
%\simeq  \fz^{\la,\reg},
%$$
The composition \eqref{comp} factors as
\begin{equation}    \label{comp1}
\fz^{\la,\nilp} \twoheadrightarrow \fz^{\la,\reg} \simeq
H^{\frac{\infty}{2}}(\fn_+\ppart,\fn_+[[t]],\BV_\la \otimes \Psi_0).
\end{equation}
\end{thm}

\begin{prop}    \label{injective}
The map
\begin{equation}    \label{inj lem}
(\BV_\la)^{\wh\n_+}_\la \to
H^{\frac{\infty}{2}}(\fn_+\ppart,\fn_+[[t]],\BV_\la \otimes \Psi_0)
\end{equation}
is injective.
\end{prop}

Assuming these two assertions, we can now prove our main result.

\medskip

\noindent {\em Proof of \thmref{weyl end}.} By \thmref{factors}, we
have the following commutative diagram:
$$
\CD
\fz^{\nilp,\lambda} @>>>  (\BV^\lambda)^{\n_+}_\la \\
@VVV                   @VVV  \\
\fz^{\reg,\lambda} @>{\sim}>>
H^{\frac{\infty}{2}}(\fn_+\ppart,\fn_+[[t]],\BV_\la \otimes \Psi_0)
\endCD
$$
The left vertical arrow is surjective, and the right vertical
arrow is injective by \propref{injective}. This readily implies that
we have an isomorphism
$$
(\BV^\lambda)^{\n_+}_\la \simeq \fz^{\reg,\lambda}.
$$
The assertion of \thmref{weyl end} follows from this, the isomorphism
\eqref{end isomo}, and the commutative diagram \eqref{important
  diagram1}.\qed

\medskip

The rest of this paper is devoted to the proof of \thmref{factors} and
\propref{injective}.

\section{Exactness}    \label{exactness}

In this section we prove \propref{injective} assuming
\thmref{factors}. The proof will rely on some properties of the
semi-infinite cohomology functors.

Let $\ghat_\crit\mod^{I,\Z}$ be the category of $\ghat_\crit$-modules
which are equivariant with respect to the Iwahori subgroup $I \subset
G\ppart$ and equipped with a $\Z$-grading with respect to the operator
$L_0 = -t\pa_t$ which commutes with $\ghat_\crit$ in the natural
way. Since $I = H \ltimes I^0$, where $I^0 = [I,I]$, the first
condition means that $\wh\n_+ = \on{Lie}(I^0)$ acts locally
nilpotently, and the constant Cartan subalgebra $\on{Lie}(H) = \h
\otimes 1 \subset \g \otimes 1 \subset \ghat_\crit$ acts semi-simply
with eigenvalues corresponding to integral weights. The second
condition means that we have an action of the extended affine
Kac--Moody algebra $\ghat'_\ka = \C L_0 \ltimes \ghat_\ka$. The
category $\ghat_\crit\mod^{I,\Z}$ is therefore the product of the
blocks of the usual category ${\mc O}_{-h^\vee}$ of modules over the
extended affine Kac--Moody algebra at the critical level $k=-h^\vee$,
corresponding to the (finite) Weyl group orbits in the set of integral
weights.

Let $\ghat_\crit\mod^{G[[t]],\Z}$ be the category of
$\ghat_\crit$-modules which are equivariant with respect to the
subgroup $G[[t]]$ and equipped with a $\Z$-grading with respect to the
operator $L_0 = -t\pa_t$. This is the full subcategory of the category
$\ghat_\crit\mod^{I,\Z}$, whose objects are modules integrable with
respect to the constant subalgebra $\g = \g \otimes 1 \subset
\ghat_\crit$.

We define $\Z$-gradings on the modules $\BM_\la$ and $\BV_\la$ in the
standard way, by setting the degrees of the generating vectors to be
equal to $0$ and using the commutation relations of $L_0$ and
$\ghat_\crit$ to define the grading on the entire modules. Thus,
$\BM_\la$ and $\BV_\la$ become objects of the category
$\ghat_\crit\mod^{I,\Z}$, and $\BV_\la$ also an object of the category
$\ghat_\crit\mod^{G[[t]],\Z}$. Moreover, the homomorphism $\BM_\la \to
\BV_\la$, and therefore the map $(\BM_\la)^{\wh\n_+}_\la \to
(\BV_\la)^{\wh\n_+}_\la$, preserve these gradings.

\medskip

We will now derive \propref{injective} from \thmref{factors} and the
following statement.

\begin{prop}    \label{exact}
The functor \eqref{n+} is right exact on the category
$\ghat_\crit\mod^{I,\Z}$ and is exact on the category
$\ghat_\crit\mod^{G[[t]],\Z}$.
\end{prop}

Introduce a $\Z$-grading operator on the standard complex of
semi-infinite cohomology
$H^{\frac{\infty}{2}+\bullet}(\fn_+\ppart,\fn_+[[t]],M \otimes
\Psi_0)$, where $M = \BM_\la$ or $\BV_\la$, by the formula
$$
L_0 - \crho \otimes 1 + \langle \la,\crho \rangle.
$$
Here $L_0$ is the natural grading operator and $\crho \in \h$ is such
that $\langle \al_\imath,\crho \rangle = 1$ for all $\imath \in I$. In
the same way as in \cite{vertex}, Sect. 15.1.8, we check that this
$\Z$-grading operator commutes with the differential of the complex
and hence induces a $\Z$-grading operator on the cohomology. Observe
that $\crho \otimes 1$ acts by multiplication by $\langle \la,\crho
\rangle$ on any element in $(\BV_\la)^{\wh\n_+}_\la$. Therefore the
map \eqref{inj lem} preserves $\Z$-gradings.

\medskip

\noindent {\em Proof of \propref{injective}.} Let
$A \in (\BV_\la)^{\wh\n_+}_\la$ be an element in the kernel of the map
\eqref{inj lem}. Since this map preserves $\Z$-gradings, without loss
of generality we may, and will, assume that $A$ is homogeneous. Under
the identification
$$
(\BV_\la)^{\wh\n_+}_\la \simeq\on{End}_{\ghat_\crit} \BV_\la
$$
it gives rise to a homogeneous $\ghat_\crit$-endomorphism $E$ of
$\BV_\la$. Then the induced map $H(E)$ on
$H^{\frac{\infty}{2}}(\fn_+\ppart,\fn_+[[t]],\BV_\la \otimes \Psi_0)$
is identically zero. Indeed, the image of the generating vector of
$(\BV_\la)^{\wh\n_+}_\la$ under the map \eqref{inj lem} is identified
with the element
$$
1 \in \fz^{\la,\reg} \simeq
H^{\frac{\infty}{2}}(\fn_+\ppart,\fn_+[[t]],\BV_\la \otimes \Psi_0),
$$
where we use the isomorphism of \thmref{factors}. Therefore the image
of this element $1$ of
$H^{\frac{\infty}{2}}(\fn_+\ppart,\fn_+[[t]],\BV_\la \otimes \Psi_0)$
under $H(E)$ is equal to the image of $A$ in
$H^{\frac{\infty}{2}}(\fn_+\ppart,\fn_+[[t]],\BV_\la \otimes \Psi_0)$,
which is $0$. By \thmref{factors},
$H^{\frac{\infty}{2}}(\fn_+\ppart,\fn_+[[t]],\BV_\la \otimes \Psi_0)$
is a free $\fz^{\la,\reg}$-module generated by the element
$1$. Therefore we find that $H(E) \equiv 0$.

Now let $M$ and $N$ be the kernel and cokernel of $E$,
$$
0 \to M \to \BV_\la \overset{E}\longrightarrow \BV_\la \to N \to 0.
$$
Note that both $M$ and $N$, as well as $\BV_\la$, are objects of the
category $\ghat_\crit\mod^{G[[t]],\Z}$. By \propref{exact}, the functor
of semi-infinite cohomology is exact on this category. Therefore we
obtain an exact sequence
\begin{multline*}
0 \to H^{\frac{\infty}{2}}(\fn_+\ppart,\fn_+[[t]],M \otimes \Psi_0)
\to H^{\frac{\infty}{2}}(\fn_+\ppart,\fn_+[[t]],\BV_\la \otimes
\Psi_0) \overset{H(E)}\longrightarrow \\
H^{\frac{\infty}{2}}(\fn_+\ppart,\fn_+[[t]],\BV_\la \otimes \Psi_0)
\to H^{\frac{\infty}{2}}(\fn_+\ppart,\fn_+[[t]],N \otimes
\Psi_0) \to 0,
\end{multline*}
where the middle map $H(E)$ is equal to zero. If $A \neq 0$, then $M$
is a proper submodule of $\BV_\la$ which does not contain the
generating vector of $\BV_\la$. We obtain that the
values of the $\Z$-grading on $M$ are strictly greater than those on
$\BV_\la$. Therefore the above sequence cannot be exact. Hence $A=0$
and we obtain the assertion of the proposition.\qed

\medskip

In the rest of this section we prove \propref{exact}. Introduce the
second semi-infinite cohomology functor (the $-$ quantum
Drinfeld--Sokolov reduction of \cite{FKW})
\begin{equation}    \label{n-}
M \mapsto H^{\frac{\infty}{2}+\bullet}(\fn_-\ppart,t\fn_-[[t]],M
\otimes \Psi_{-\crho}),
\end{equation}
where
\begin{equation}    \label{Psirho}
\Psi^-_{-\crho}: \n_-\ppart \to \C
\end{equation}
is given by the formula
$$
\Psi^-_{-\crho}(f_{\al,n}) = \begin{cases} 1, & \on{if} \al =
  \al_\imath, n=0, \\ 0, & \on{otherwise}. \end{cases}
$$
Then we have the following important result due to Arakawa \cite{Ar},
Main Theorem 1,(1) (note that the functor \eqref{n-} is
the functor $H^\bullet_-$ in the notation \cite{Ar}):

\begin{thm}    \label{arakawa}
The functor \eqref{n-} is exact on the category
$\ghat_\crit\mod^{I,\Z}$.
\end{thm}

We now derive \propref{exact} from \thmref{arakawa}.

\medskip

\noindent {\em Proof of \propref{exact}.} Recall that we have the
convolution functors
$$
M \mapsto {\mc F} \star M
$$
on the category of $I$-equivariant $\ghat_\crit$ modules, for each
$I$-equivariant right D-module ${\mc F}$ on $G\ppart/I$ (see
\cite{FG:local}, Sect. 22, for the precise definition).

According to Proposition 18.1.1 of \cite{FG:local}, we have
\begin{equation}    \label{isom n- and n+}
H^{\frac{\infty}{2}+\bullet}(\fn_+\ppart,\fn_+[[t]],M
\otimes \Psi_0) \simeq
H^{\frac{\infty}{2}+\bullet}(\fn_-\ppart,t\fn_-[[t]],j_{w_0\crho,*}
\star M \otimes \Psi_{-\crho})
\end{equation}
for any $I$-equivariant $\ghat_\crit$-module $M$. We recall that the
D-module $j_{w_0\crho,*}$ is defined as the $*$-extension of the
``constant'' D-module on the $I$-orbit in the affine flag scheme
$G\ppart/I$ corresponding to the element $w_0\crho$ of the affine
Weyl group. Hence the functor
$$
M \mapsto j_{w_0\crho,*}\star M
$$
is right exact. Combining this with \thmref{arakawa}, we obtain that
the functor
$$
M \mapsto
H^{\frac{\infty}{2}+\bullet}(\fn_-\ppart,t\fn_-[[t]],j_{w_0\crho,*}
\star M \otimes \Psi_{-\crho})
$$
is right exact on the category $\ghat_\crit\mod^{I,\Z}$ (note that the
convolution with $j_{w_0\crho,*}$ sends $\Z$-graded modules to
$\Z$-graded modules). The isomorphism \eqref{isom n- and n+} then
implies that the functor \eqref{n+} is right exact on the category
$\ghat_\crit\mod^{I,\Z}$.

On the other hand, let $j_{w_0\crho,!}$ be the $!$-extension of
the ``constant'' D-module on the same $I$-orbit. The convolution
functor with $j_{w_0\crho,!}$ is both left and right adjoint to
the convolution with $j_{w_0\crho,*}$. Therefore we find that the
functor
$$
M \mapsto j_{w_0\crho,!}\star M
$$
is left exact. Combining this with \thmref{arakawa}, we obtain that
the functor
$$
M \mapsto
H^{\frac{\infty}{2}+\bullet}(\fn_-\ppart,t\fn_-[[t]],j_{w_0\crho,!}
\star M \otimes \Psi_{-\crho})
$$
is left exact on the category $\ghat_\crit\mod^{I,\Z}$ (again, note
that the convolution with $j_{w_0\crho,!}$ sends $\Z$-graded
modules to $\Z$-graded modules).

Now consider the homomorphism
\begin{equation}    \label{isomor}
j_{w_0\crho,!}\star M \to j_{w_0\crho,*}\star M
\end{equation}
of $\ghat_\crit$-modules induced by the morphism
$$
j_{w_0\crho,!} \to j_{w_0\crho,*}
$$
of D-modules on $G\ppart/I$. Suppose in addition that $M$ is
$G[[t]]$-equivariant. Then we have
$$
j_{w_0\crho,!}\star M \simeq j_{w_0\cdot \check\rho,!} \star
\delta_{1_{Gr_G}} \underset{G[[t]]}\star M, \qquad j_{w_0\crho,*}
\star M \simeq j_{w_0\cdot \check\rho,*} \star \delta_{1_{Gr_G}}
\underset{G[[t]]}\star M,
$$
where $\delta_{1_{Gr_G}}$ is the ``delta-function'' D-module on the
affine Grassmannian $Gr_G = G\ppart/G[[t]]$ supported at the identity
coset. It then follows from Lemma 15.1.2 of \cite{FG:local} that the
kernel and the cokernel of the map \eqref{isomor} are partially
integrable $\ghat_\crit$-modules. We recall from \cite{FG:local},
Sect. 6.3, that a $\ghat_\crit$-module is called partially integrable
if it admits a filtration such that each successive
quotient is equivariant with respect to the parahoric Lie
subalgebra $\fp^\iota=\Lie(I)+{\mathfrak {sl}}_2^\iota$ for some
vertex of the Dynkin diagram of $\g$, $\iota\in I$.

But, according to Lemma 18.1.2 of \cite{FG:local}, if $M$ is a
partially integrable $\ghat_\crit$-module, then
$$
H^{\frac{\infty}{2}+\bullet}(\fn_-\ppart,t\fn_-[[t]],M \otimes
\Psi_{-\crho}) = 0
$$
for all $i \in \Z$. Therefore we obtain that the map
$$
H^{\frac{\infty}{2}+\bullet}(\fn_-\ppart,t\fn_-[[t]],j_{w_0\crho,!}
\star M \otimes \Psi_{-\crho}) \to
H^{\frac{\infty}{2}+\bullet}(\fn_-\ppart,t\fn_-[[t]],j_{w_0\crho,*}
\star M \otimes \Psi_{-\crho})
$$
induced by \eqref{isomor} is an isomorphism for any
$G[[t]]$-equivariant $\ghat_\crit$-module $M$. Since the former is
right exact and the latter is left exact on the category
$\ghat_\crit\mod^{G[[t]],\Z}$, we obtain that both functors are exact
on this category. Combining this with the isomorphism \eqref{isom n-
and n+}, we find that the functor
$$
M \mapsto H^{\frac{\infty}{2}+\bullet}(\fn_+\ppart,\fn_+[[t]],M
\otimes \Psi_0)
$$
is exact on the category $\ghat_\crit\mod^{G[[t]],\Z}$.

This completes the proof of \propref{exact}.\qed

\medskip

\begin{remark} There are obvious analogues of the categories
$\ghat_\crit\mod^{I,\Z}$ and $\ghat_\crit\mod^{G[[t]],\Z}$ for an
arbitrary level $\ka$. The same proof as above works for any $\ka$, so
\propref{exact} actually holds for an arbitrary level.\qed
\end{remark}

\medskip

It remains to prove \thmref{factors}. We will give two proofs: one
relies on the results of \cite{FG:fusion}, and the other uses the
Wakimoto modules. Both proofs use the computation of the characters of
$\fz^{\la,\reg}$ and
$H^{\frac{\infty}{2}+\bullet}(\fn_+\ppart,\fn_+[[t]],\BV_\la \otimes
\Psi_0)$ which is performed in the next section.

\section{Computation of characters}    \label{char z}

\subsection{Character of $\fz^{\la,\reg}$}

Let us compute the character of the algebra
$$
\fz^{\la,\reg} = \on{Fun} \on{Op}_{\cg}^{\la,\reg}
$$
of functions on $\on{Op}_{\cg}^{\la,\reg}$ with respect to the
$\Z$-grading by the operator $L_0$. This space was defined in
\secref{first section}, but note that we now switch to the Langlands
dual Lie algebra $\cg$.

We give another, more convenient, realization of the space
$\on{Op}_{\g}^{\la,\reg}$. Suppose that we are given an operator of
the form
\begin{equation}    \label{oper on disc reg}
\nabla = \pa_t + \sum_{\imath \in I} t^{\langle \chal_\imath,\la
\rangle} f_\imath + \bv(t), \qquad \bv(t)
\in \check\fb(\hCO).
\end{equation}
Applying gauge transformation with $(\la+\rho)(t)$, we obtain an
operator of the form
\begin{equation}    \label{oper on disc reg1}
\nabla' = \pa_t + \frac{1}{t} \left( \sum_{\imath \in I} f_\imath -
(\la+\rho) \right) + \bv(t), \qquad \bv(t) \in
(\la+\rho)(t)\check\fb(\hCO)(\la+\rho)(t)^{-1}.
\end{equation}
The space $\on{Op}_{\cg}^{\la,\reg}$ is defined as the space of
$\check{N}(\hCO)$-equivalence classes of operators \eqref{oper on disc
reg}. Equivalently, this is the space of $(\la+\rho)(t)
\check{N}(\hCO) (\la+\rho)(t)^{-1}$-equivalence classes of operators
\eqref{oper on disc reg1}. It follows from Theorem 2.21 of
\cite{FG:local} that the action of the group $(\la+\rho)(t)
\check{N}(\hCO) (\la+\rho)(t)^{-1}$ on this space is free. Therefore
the character of $\on{Fun} \on{Op}_{\cg}^{\la,\reg}$ is equal to the
the character of the algebra of functions on the space of operators
\eqref{oper on disc reg1} divided by the character of the algebra
$\on{Fun}((\la+\rho)(t) \check{N}(\hCO) (\la+\rho)(t)^{-1})$.

By definition, a character of a $\Z_+$-graded vector space $V =
\bigoplus_{n \in \Z_+} V_n$, where each $V_n$ is finite-dimensional,
is the formal power series
$$
\on{ch} V = \sum_{n\geq 0} \dim V_n \cdot q^n.
$$

Applying the dilation $t \mapsto at$ to \eqref{oper on disc reg1}, we
obtain the action $\bv(t) \mapsto a\bv(at)$. Therefore
the character of the algebra of functions on the space of
operators \eqref{oper on disc reg1} is equal to
$$
\prod_{n>0} (1-q^n)^{-\ell} \cdot \prod_{\chal \in \check\De_+}
\prod_{n>0} (1-q^{n+\langle \chal,\la+\rho \rangle})^{-1}.
$$
On the other hand, the character of $\on{Fun}((\la+\rho)(t)
\check{N}(\hCO) (\la+\rho)(t)^{-1})$ is equal to
$$
\prod_{\chal \in \check\De_+}
\prod_{n\geq 0} (1-q^{n+\langle \chal,\la+\rho \rangle})^{-1}.
$$
Therefore the character of $\fz^{\la,\reg} = \on{Fun}
\on{Op}_{\cg}^{\la,\reg}$ is equal to
$$
\prod_{n>0} (1-q^n)^{-\ell} \cdot \prod_{\chal \in \check\De_+}
(1-q^{\langle \al,\la+\rho \rangle}).
$$

We rewrite this in the form
\begin{equation}    \label{char Op la}
\on{ch} \fz^{\la,\reg} = \prod_{\chal \in \check\Delta_+}
\frac{1-q^{\langle \chal,\la+\rho \rangle}}{1-q^{\langle \chal,\rho
\rangle}} \prod_{i=1}^\ell \prod_{n_i \geq d_i+1} (1-q^{n_i})^{-1},
\end{equation}
using the identify
$$
\prod_{\chal \in \check\De_+} (1-q^{\langle \chal,\rho \rangle}) =
\prod_{i=1}^\ell \prod_{m_i=1}^{d_i} (1-q^{m_i}),
$$
where $d_1,\ldots,d_\ell$ are the exponents of $\g$.

\subsection{Computation of semi-infinite cohomology}

Let us now compute the semi-infinite cohomology of $\BV_\la$. The
complex $C^{\bullet}(\BV_\la)$ computing this cohomology is described
in \cite{vertex}, Ch. 15. In particular, as explained in
\secref{exactness}, it carries a $\Z$-grading operator which commutes
with the differential and hence gives rise to a grading operator on
the cohomology. We will compute the character with respect to this
grading operator.

\begin{thm}    \label{si coh}
We have
$$
\on{ch} H^{\frac{\infty}{2}}(\fn_+\ppart,\fn_+[[t]],\BV_\la \otimes
\Psi_0) = \on{ch} \fz^{\la,\reg}
$$
given by formula \eqref{char Op la}, and
$$
H^{\frac{\infty}{2}+i}(\fn_+\ppart,\fn_+[[t]],\BV_\la \otimes
\Psi_0) = 0, \qquad i \neq 0.
$$
\end{thm}

\begin{proof}
The vanishing the $i$th cohomology for $i \neq 0$ follows from
\propref{exact}. We will give an alternative proof of this, as well as
the computation of the character of the $0$th cohomology, using the
argument of \cite{vertex}, Sect. 15.2.

Consider the complex $C^{\bullet}(\BV_\la)$ computing our
semi-infinite cohomology. This complex was studied in detail in
\cite{vertex}, Sect. 15.2, in the case when $\la=0$. We decompose
$C^{\bullet}(\BV_\la)$ into the tensor product of two subcomplexes as
in \cite{vertex}, Sect. 15.2.1,
$$
C^{\bullet}(\BV_\la) = C^{\bullet}(\BV_\la)_0 \otimes
C^{\bullet}(\BV_\la)'
$$
defined in the same way as the subcomplexes
$C^\bullet_{-h^\vee}(\g)_0$ and $C^\bullet_{-h^\vee}(\g)'$,
respectively. In fact, $C^{\bullet}(\BV_\la)_0$ is equal to
$C^\bullet_{-h^\vee}(\g)_0$, and
$$
C^{\bullet}(\BV_\la)' \simeq V_\la \otimes U(t^{-1}\fb_-[t^{-1}])
\otimes \bigwedge{}^\bullet(\fn_+[[t]]^*).
$$
In particular, its cohomological grading takes only non-negative
values on $C^{\bullet}(\BV_\la)'$.

We show, in the same way as in \cite{vertex}, Lemma 15.2.5, that the
cohomology of $C^{\bullet}(\BV_\la)$ is isomorphic to the tensor
product of the cohomologies of the subcomplexes
$C^{\bullet}(\BV_\la)_0$ and $C^{\bullet}(\BV_\la)'$. The former is
one-dimensional, according to \cite{vertex}, Lemma 15.2.7, and hence
we find that our semi-infinite cohomology is isomorphic to the
cohomology of the subcomplex $C^{\bullet}(\BV_\la)'$.

Following verbatim the computation in
\cite{vertex}, Sect. 15.2.9, in the case when $\la=0$, we find that
the $0$th cohomology of $C^{\bullet}(\BV_\la)'$ is isomorphic to
$$
H^{\frac{\infty}{2}}(\fn_+\ppart,\fn_+[[t]],\BV_\la \otimes
\Psi_0) \simeq V_\la \otimes V(\fa_-)
$$
(where $V(\fa_-)$ is defined in \cite{vertex}, Sect. 15.2.9), and all
other cohomologies vanish.

In particular, we find that the character of
$H^{\frac{\infty}{2}}(\fn_+\ppart,\fn_+[[t]],\BV_\la \otimes
\Psi_0)$ is equal to
$$
\on{ch} V_\la \cdot \on{ch} V(\fa_-),
$$
where $\on{ch} V_\la$ is the character of $V_\la$ with respect to the
principal grading. By \cite{vertex}, Sect. 15.2.9, we have
$$
\on{ch} V(\fa_-) = \prod_{i=1}^\ell
\prod_{m_i \geq d_i+1} (1-q^{m_i})^{-1}.
$$
According to formula (10.9.4) of \cite{Kac}, $\on{ch} V_\la$ is equal
to
\begin{equation}    \label{Vla char2}
\on{ch} V_\la = \prod_{\check\al \in \check\Delta_+}
\frac{1-q^{\langle \chal,\la+\rho \rangle}}{1-q^{\langle \chal,\rho
    \rangle}}.
\end{equation}
Therefore the character
$$
\on{ch} H^{\frac{\infty}{2}}(\fn_+\ppart,\fn_+[[t]],\BV_\la \otimes
\Psi_0)
$$
is given by formula \eqref{char Op la}, which coincides with the
character of $\fz^{\la,\reg}$.
\end{proof}

\section{Proof of \thmref{factors}}    \label{proof}

\subsection{First proof}    \label{first proof}

The following result is proved in \cite{FG:fusion}, Lemma 1.7.

\begin{prop}    \label{is contained}
The action of the center $\fZ_{\fg}$ on $\BV_\la$ factors through
$\fz^{\la,\reg}$.
\end{prop}

Let $I_\la$ be the ideal of $\Op^{\reg,\la} = \on{Spec}
\fz^{\reg,\la}$ in the center $\fZ_{\g} =
\on{Fun}(\Op(\D^\times))$. As explained in \cite{FG:local},
Sect. 4.6 (see \cite{BD}, Sect. 3.6, in the case when $\la=0$), the
Poisson structure on $\fZ_{\g}$ gives rise to the structure of a Lie
algebroid on the quotient $I_\la/(I_\la)^2$, which we denote by
$N^*_{\Op^{\la,\reg}/\on{Op}_{\cg}(\D^\times)}$. Recall from
\cite{FF,F:wak} that the Poisson structure on $\fZ_{\g}$ is
obtained by deforming the completed enveloping algebra of $\ghat$ to
non-critical levels. By \propref{is contained}, $I_\la$ annihilates
the module $\BV_\la$. Since this module may be deformed to the Weyl
modules at non-critical levels, we obtain that the Lie algebroid
$N^*_{\Op^{\la,\reg}/\on{Op}_{\cg}(\D^\times)}$ naturally acts on
$\BV_\la$ (see \cite{BD}, Sect. 5.6, in the case when $\la=0$) and on
its semi-infinite cohomology. This action is compatible with the
action of $\on{Fun}(\Op(\D^\times))/I_\la = \fz^{\reg,\la}$.

Using the commutative diagram \eqref{important diagram1} and
\propref{is contained}, we obtain that the composition \eqref{comp}
factors through a map
$$
\fz^{\la,\reg} \to
H^{\frac{\infty}{2}}(\fn_+\ppart,\fn_+[[t]],\BV_\la \otimes \Psi_0).
$$
By applying the same argument as in the proof of Proposition 18.3.2 of
\cite{FG:local}, we obtain that the above map is a homomorphism of
modules over the Lie algebroid
$N^*_{\Op^{\la,\reg}/\on{Op}_{\cg}(\D^\times)}$. This homomorphism is
clearly non-zero, because we can identify the image of the generator
$1 \in \fz^{\la,\reg}$ with the cohomology class represented by the
highest weight vector of $\BV_\la$. Since $\fz^{\reg,\lambda}$ is
irreducible as a module over
$N^*_{\Op^{\la,\reg}/\on{Op}_{\cg}(\D^\times)}$, this homomorphism is
injective. Moreover, it is clear that this map preserves the natural
$\Z$-gradings on both modules. Therefore the equality of the two
characters established in \secref{char z} shows that it is an
isomorphism.\qed

\subsection{Second proof}

Let us recall some results about Wakimoto modules of critical level
from \cite{F:wak} (see also \cite{FG:local}). Specifically, we will
consider the module which was denoted by $W_{\la,\ka_c}$ in
\cite{F:wak} and by $\BW^{w_0}_{\crit,\la}$ in \cite{FG:local}. Here
we will denote it simply by $W_\la$. As a vector space, it is
isomorphic to the tensor product
\begin{equation}    \label{mod Wla1}
W_\la = M_\g \otimes \fH^\la,
\end{equation}
where we use the notation
$$
\fH^\la = \on{Fun} \ConHD^{\RS,-\la}
$$
(see \secref{first section}), and $M_\g$ is the Fock representation
of a Weyl algebra with generators $a_{\al,n}, a^*_{\al,n}, \al \in
\De_+, n \in \Z$.

We now construct a map $\BV_\la \to W_{\la}$. Let us observe that the
action of the constant subalgebra $\g \subset \ghat_\ka$ on the
subspace
$$
W_\la^0 = \C[a^*_{\al,0}]_{\al \in \De_+} \vac \simeq \on{Fun} N
\subset M_\g
$$
coincides with the natural action of $\g$ on the contragredient Verma
module $M^*_\la$ realized as $\on{Fun} N$. In addition,
the Lie subalgebra $\g \otimes t\otimes \C[[t]] \subset \ghat_\crit$
acts by zero on the subspace $W^0_\la$, and ${\bf 1}$ acts as the
identity.

Therefore the injective $\g$-homomorphism $V_\la \hookrightarrow
M^*_\la$ gives rise to a non-zero $\ghat_\crit$-homomorphism
$$
\imath_\la: \BV_\la \to W_\la
$$
sending the generating subspace $V_\la \subset \BV_\la$ to the image
of $V_\la$ in $W_\la^0 \simeq M^*_\la$.

Now we obtain a sequence of maps
\begin{equation}    \label{seq of maps}
\BM_{\la} \twoheadrightarrow \BV_\la \rightarrow W_{\la}.
\end{equation}
Recall from \secref{hom ghat} that
$$
(\BM_\la)^{\wh\n_+}_\la \simeq \fz^{\la,\nilp}.
$$
We also prove, by using the argument of Lemma 6.5 of \cite{F:wak} that
$$
(W_\la)^{\wh\n_+}_\la = \fH^\la,
$$
where $\fH^\la$ is identified with the second factor of the
decomposition \eqref{mod Wla1}.

We obtain the following commutative diagram, in which all maps preserve
$\Z$-gradings:
\begin{equation}    \label{one more}
\begin{CD}
\fz^{\la,\nilp} @>{b}>> \fH^\la  \\ @V{a}VV @V{c}VV \\
H^{\frac{\infty}{2}}(\fn_+\ppart,\fn_+[[t]],\BV_\la \otimes
\Psi_0) @>{d}>> H^{\frac{\infty}{2}}(\fn_+\ppart,\fn_+[[t]],W_\la
\otimes \Psi_0)
\end{CD}
\end{equation}
Here the map $a$ is obtained as the composition
$$
\fz^{\la,\nilp} \simeq (\BM_\la)^{\wh\n_+}_\la \to
(\BV_\la)^{\wh\n_+}_\la \to
H^{\frac{\infty}{2}}(\fn_+\ppart,\fn_+[[t]],\BV_\la \otimes \Psi_0),
$$
the map $c$ as the composition
$$
\fH^\la \simeq (W_\la)^{\wh\n_+}_\la \to
H^{\frac{\infty}{2}}(\fn_+\ppart,\fn_+[[t]],W_\la \otimes \Psi_0)
$$
(see \lemref{map}), and the map $b$ as the composition
$$
\fz^{\la,\nilp} \simeq (\BM_\la)^{\wh\n_+}_\la \to
(W_\la)^{\wh\n_+}_\la \simeq \fH^\la.
$$

Using Theorem 12.5 of \cite{F:wak}, we identify the map $b$ with the
homomorphism of the algebras of functions corresponding to the map
$$
\ConHD^{\RS,-\la} \to \on{Op}_{\cg}^{\la,\nilp}
$$
obtained by restriction from the Miura transformation
\eqref{MT}. Therefore we obtain from \lemref{prin N bdle} that the
image of this map is precisely $\on{Op}_{\cg}^{\la,\reg}$, and so the
image of the homomorphism $b$ is equal to
$$
\fz^{\la,\reg} = \on{Fun} \on{Op}_{\cg}^{\la,\reg}.
$$

On the other hand, we have the decomposition \eqref{mod Wla1} and
$\n_+\ppart$ acts along the first factor $M_\g$. It follows from the
definition of $M_\g$ that there is a canonical identification
$$
H^{\frac{\infty}{2}}(\fn_+\ppart,\fn_+[[t]],M_\g \otimes
\Psi_0) \simeq \C \vac.
$$
Therefore we obtain that
$$
H^{\frac{\infty}{2}}(\fn_+\ppart,\fn_+[[t]],W_\la
\otimes \Psi_0) \simeq \fH^\la
$$
and so the map $c$ is an isomorphism.

This implies that the image of the composition $c \circ b$ is
$\fz^{\la,\reg} \subset \fH^\la$. Therefore $d$ factors as follows:
$$
H^{\frac{\infty}{2}}(\fn_+\ppart,\fn_+[[t]],\BV_\la
\otimes \Psi_0) \twoheadrightarrow \fz^{\la,\reg} \hookrightarrow
\fH^\la.
$$
But the characters of the first two spaces coincide, according to
\thmref{si coh}. Hence
$$
H^{\frac{\infty}{2}}(\fn_+\ppart,\fn_+[[t]],\BV_\la
\otimes \Psi_0) \simeq \fz^{\la,\reg}.
$$
This proves \thmref{factors}.\qed

\bigskip

We obtain the following corollary, which for $\la=0$ was proved in
\cite{F:wak}, Prop. 5.2.

\begin{cor}    \label{imath}
The map $\imath_\la: \BV_\la \to W_\la$ is injective for any dominant
integral weight $\la$.
\end{cor}

\begin{proof}
Let us extend the action of $\ghat_\ka$ to an action of $\ghat'_\ka =
\C L_0 \ltimes \ghat_\ka$ in the same way as above. Denote by $K_\la$
the kernel of the map $\imath_\la$. Suppose that $K_\la \neq 0$. Since
both $\BV_\la$ and $\BM_\la$ become graded with respect to the
extended Cartan subalgebra $\wt\h = (\h \otimes 1) \oplus \C L_0$, we
find that the $\ghat_\crit$-module $K_\la$ contains a non-zero highest
weight vector annihilated by the Lie subalgebra $\wh\n_+$. Let
$\wh\mu$ be its weight. Since $\BV_\la$ is $\g$-integrable, so is
$K_\la$, and therefore the restriction of $\wt\mu$ to $\h \subset
\wt\h$ is a dominant integral weight $\mu$.

Now, since $\BV_\la$ is a quotient of the Verma module $\BM_\la$, and
the action of the center $\fZ_{\g}$ on $\BM_\la$ factors through
$\on{Fun} \Op^{\la,\nilp}$ (see the diagram \eqref{important
diagram1}), we find that the same is true for $\BV_\la$.\footnote{We
have already determined in \propref{is contained} that the support of
$\BV_\la$ is contained in $\Op^{\la,\reg}$, but it is not necessary to
use this result here.} According to \cite{F:wak}, Prop. 12.8, the
degree $0$ part of $\on{Fun} \Op^{\la,\nilp}$ is isomorphic to
$(\on{Fun} \h^*)^W$, and it acts on any highest weight vector $v$
through the quotient by the maximal ideal in $(\on{Fun} \h^*)^W$
corresponding to $-\la-\rho$. Moreover, its action coincides with the
action of the center $Z(\g)$ of $U(\g)$ on the vector $v$ under the
Harish-Chandra isomorphism $Z(\g) \simeq (\on{Fun} \h^*)^W$ followed
by the sign isomorphism $x \mapsto -x$. Therefore if $v$ has
$\h$-weight $\mu$, the action of $(\on{Fun} \h^*)^W$ factors through
the quotient by the maximal ideal in $(\on{Fun} \h^*)^W$ corresponding
to $-\mu-\rho$. Since $\mu$ is dominant integral, this implies that
$\mu=\la$.

Thus, $K_\la$ contains a non-zero highest weight vector of $\h$-weight
$\la$. This vector then lies in $(\BV_\la)^{\wh\n_+}_\la$. But
\thmref{weyl end} and the diagram \eqref{one more} imply that the map
$$
(\BV_\la)^{\wh\n_+}_\la \to W_\la
$$
is injective. Indeed, the image of the composition
$$
(\BM_\la)^{\wh\n_+}_\la \to (\BV_\la)^{\wh\n_+}_\la \to W_\la
$$
is equal to $\fz^{\la,\reg} \subset \fH^\la$, which is isomorphic to
$(\BV_\la)^{\wh\n_+}_\la$.

This leads to a contradiction and hence proves the
desired assertion.
\end{proof}

The following result is also useful in applications. Note that for
$\la=0$ it follows from the corresponding statement for the associated
graded module proved in \cite{EF} or from the results of \cite{BD},
Sect. 6.2.

\begin{thm}    \label{free}
For any dominant integral weight $\la$ the Weyl module $\BV_\la$ is
free as a $\on{Fun} \Op^{\la,\reg}$-module.
\end{thm}

\begin{proof}
Recall from \secref{first proof} that $\BV_\la$ carries an action of
the Lie algebroid $N^*_{\Op^{\la,\reg}/\on{Op}_{\cg}(\D^\times)}$
compatible with the action of $\fz^{\reg,\la} =
\on{Fun}(\Op^{\la,\reg})$. According to the results of
\cite{FG:local}, Sect. 4.6, this Lie algebroid is nothing but the
Atiyah algebroid of the universal $\cG$-bundle on $\Op^{\la,\reg}$
whose fiber at $\chi \in \Op^{\la,\reg}$ is the fiber of the
$\cG$-torsor on $\D$ underlying $\chi$ at $0 \in \D$. This bundle is
isomorphic (non-canonically) to the trivial $\cG$-bundle on
$\Op^{\la,\reg}$. Let us choose such an isomorphism. Then
$N^*_{\Op^{\la,\reg}/\on{Op}_{\cg}(\D^\times)}$ splits as a direct sum
of the Lie algebra $\on{Vect}(\Op^{\la,\reg})$ of vector fields on
$\Op^{\la,\reg}$ and the Lie algebra $\cg \otimes
\on{Fun}(\Op^{\la,\reg})$.

Thus, we obtain an action of $\on{Vect}(\Op^{\la,\reg})$ on $\BV_\la$
compatible with the action of $\on{Fun}(\Op^{\la,\reg})$. Note that
the algebra $\fz^{\reg,\la}$ and the Lie algebra
$\on{Vect}(\Op^{\la,\reg})$ are $\Z$-graded by the operator $L_0 = -
t\pa_t$. According to \lemref{prin N bdle}, $\Op^{\la,\reg}$ is an
infinite-dimensional affine space, and there exists a system of
coordinates $x_i, i=1,2,\ldots$, on it such that these coordinates are
homogeneous with respect to $L_0$. The character formula \eqref{char
Op la} shows that the degrees of all of these generators are strictly
positive.

Thus, we find that the action on $\BV_\la$ of the polynomial algebra
$\on{Fun}(\Op^{\la,\reg})$ generated by the $x_i$'s extends to an
action of the Weyl algebra generated by the $x_i$'s and the $\pa/\pa
x_i$'s. Recall that the $\Z$-grading on $\BV_\la$ with respect to the
operator $L_0$ takes non-negative integer values. Repeating the
argument of Lemma 6.2.2 of \cite{BD} (see Lemma 9.13 of \cite{Kac}),
we obtain that $\BV_\la$ is a free module over $\fz^{\reg,\la} =
\on{Fun}(\Op^{\la,\reg})$.
\end{proof}

\end{document}